\documentclass[11.05pt,a4paper]{amsart}
\usepackage{amssymb,amsmath}
\usepackage{color}
\usepackage{latexsym}
\usepackage{amsthm,amsfonts,amssymb,mathrsfs}
\usepackage{rotating}
\usepackage[leqno]{amsmath}
\usepackage{xspace}
\usepackage[all]{xy}
\usepackage{longtable}

\headsep=1cm \numberwithin{equation}{section}
\newtheorem{theorem}{Theorem}[section]
\newtheorem{lemma}[theorem]{Lemma}
\newtheorem{proposition}[theorem]{Proposition}
\newtheorem{corollary}[theorem]{Corollary}

\theoremstyle{definition}

\theoremstyle{remark}
\newtheorem{remark}[theorem]{Remark}

\newtheorem{acknowledgement}{Acknowledgement}

\newcommand{\Ass}{\operatorname{Ass}}

\newcommand{\Assh}{\operatorname{Assh}}
\newcommand{\Spec}{\operatorname{Spec}}

\newcommand{\amp}{\operatorname{amp}}

\newcommand{\Ht}{\operatorname{ht}}

\newcommand{\V}{\operatorname{V}}

\newcommand{\Ext}{\operatorname{Ext}}
\newcommand{\Supp}{\operatorname{Supp}}

\newcommand{\Hom}{\operatorname{Hom}}
\newcommand{\Att}{\operatorname{Att}}
\newcommand{\Ann}{\operatorname{Ann}}

\newcommand{\Vdim}{\operatorname{Vdim}}
\newcommand{\depth}{\operatorname{depth}}

\newcommand{\Max}{\operatorname{Max}}

\newcommand{\lo}{\longrightarrow}
\newcommand{\fm}{\frak{m}}
\newcommand{\fp}{\frak{p}}
\newcommand{\fq}{\frak{q}}
\newcommand{\fa}{\frak{a}}
\newcommand{\fb}{\frak{b}}

\newcommand{\fn}{\frak{n}}

\newenvironment{prf}[1][Proof]{\begin{proof}[\bf #1]}{\end{proof}}
\begin{document}

\author[M. Hatamkhani and K. Divaani-Aazar]{Marziyeh Hatamkhani and Kamran Divaani-Aazar}
\title[The derived category analogue of the ...]
{The derived category analogue of the Hartshorne-Lichtenbaum
Vanishing Theorem}

\address{M. Hatamkhani, Department of Mathematics, Az-Zahra University,
Vanak, Post Code 19834, Tehran, Iran} \email{hatamkhanim@yahoo.com}

\address{K. Divaani-Aazar, Department of Mathematics, Az-Zahra University,
Vanak, Post Code 19834, Tehran, Iran-and-School of Mathematics,
Institute for Research in Fundamental Sciences (IPM), P.O. Box
19395-5746, Tehran, Iran} \email{kdivaani@ipm.ir}


\keywords {Attached prime ideals; derived categories; local
cohomology; stable under specialization subsets.\\
The second author was supported by a grant from IPM (No. 90130212).}

\begin{abstract} Let $\fa$ be an ideal of a local ring $(R,\fm)$
and $X$ a $d$-dimensional homologically bounded complex of
$R$-modules whose all homology modules are finitely generated. We
show that $H^d_{\fa}(X)=0$ if and only if $\dim \widehat{R}/\fa
\widehat{R}+\fp>0$ for all prime ideals $\fp$ of $\hat{R}$ such that
$\dim \hat{R}/\fp-\inf (X\otimes_R\hat{R})_{\fp}=d$.
\end{abstract}

\maketitle

\section{Introduction}

The Hartshorne-Lichtenbaum Vanishing Theorem is one of the most
important results in the theory of local cohomology modules. There
are several proofs known now of this result; see e.g. \cite{BH},
\cite{CS} and \cite{Sc}. Also, there are several generalizations of
this result. The second named author, Naghipour and Tousi \cite{DNT}
have extended it to local cohomology with support in stable under
specialization subsets.  Takahashi, Yoshino and Yoshizawa \cite{TYY}
have extended it to local cohomology with respect to pairs of
ideals. Also, more recently, the Hartshorne-Lichtenbaum Vanishing
Theorem is extended to generalized local cohomology modules; see
\cite{DH}. Our aim in this paper is to establish a generalization
of the Hartshorne-Lichtenbaum Vanishing Theorem which contains all
of these generalizations. We do this by establishing the derived
category analogue of the Hartshorne-Lichtenbaum Vanishing Theorem.
For giving the precise statement of this result, we need to fix some
notation.

Throughout, $R$ is a commutative Noetherian ring with nonzero identity.
The derived category of $R$-modules is denoted by
$\mathcal{D}(R)$. We use the symbol $\simeq$ for denoting
isomorphisms in $\mathcal{D}(R)$. For a complex $X\in
\mathcal{D}(R)$, its supremum and infimum are defined, respectively,
by $\sup X:=\sup \{i\in \mathbb{Z}|H_i(X)\neq 0\}$ and $\inf X:=\inf
\{i\in \mathbb{Z}|H_i(X) \neq 0\}$, with the usual convention that
$\sup \emptyset=-\infty$ and $\inf \emptyset=\infty$. Also,
amplitude of $X$ is defined by $\amp X:=\sup X-\inf X$. Recall that
$\dim_RX$ is defined by $\dim_RX:=\sup\{\dim R/\fp-\inf
X_{\fp}|\fp\in \Spec R\}$ and we define $\Assh_RX$ by
$$\Assh_RX:=\{\fp\in \Spec R|\dim R/\fp-\inf X_{\fp}=\dim_RX \}.$$

Any $R$-module $M$ can be considered as a complex having $M$ in its 0-th
spot and 0 in its other spots. We denote the full subcategory of homologically
left bounded complexes by $\mathcal{D}_{\sqsubset}(R)$. Also, we denote the full
subcategory of complexes with finitely generated homology modules that are
homologically bounded (resp. homologically left bounded) by $\mathcal{D}_{\Box}^f(R)$
(resp. $\mathcal{D}_{\sqsubset}^f(R)$).

Let $\fa$ be an ideal of $R$ and $X\in \mathcal{D}_{\sqsubset}(R)$. A subset $\mathcal{Z}$ of
$\Spec R$ is said to be {\it stable under specialization} if $\V(\fp)\subseteq \mathcal{Z}$
for all $\fp \in \mathcal{Z}$. For any $R$-module $M$, $\Gamma_{\mathcal{Z}}(M)$ is defined by
$$\Gamma_{\mathcal{Z}}(M):=\{x\in M|\Supp_RRx\subseteq \mathcal{Z}\}.$$ The right derived
functor of the functor $\Gamma_{\mathcal{Z}}(-)$ exists in $\mathcal{D}(R)$ and the complex ${\bf
R}\Gamma_{\mathcal{Z}}(X)$ is defined by ${\bf R}
\Gamma_{\mathcal{Z}}(X):=\Gamma_{\mathcal{Z}}(I)$, where $I$ is any
injective resolution of $X$.  Also, for any
integer $i$, the $i$-th local cohomology module of $X$ with respect
to $\mathcal{Z}$ is defined by $H_{\mathcal{Z}}^i(X):=H_{-i}({\bf R}
\Gamma_{\mathcal{Z}}(X))$.  To comply with the usual notation, for
$\mathcal{Z}:=\V(\fa)$, we denote ${\bf R}\Gamma_{\mathcal{Z}}(-)$
and $H_{\mathcal{Z}}^i(-)$ by ${\bf R}\Gamma_{{\fa}}(-)$ and
$H_{\fa}^i(-)$, respectively. By \cite[Corollary 3.7 and Proposition 3.14 d)]{F3},
for any complex $X\in \mathcal{D}_{\Box}^f(R)$, we know
that $$\sup \{i\in \mathbb{Z}|H_{\fa}^i(X)\neq 0 \}\leq \dim_RX$$
with equality if $R$ is local and $\fa$ is its maximal ideal. Denote
the set of all ideals $\fb$ of $R$ such that $\V(\fb)\subseteq
\mathcal{Z}$ by $F(\mathcal{Z})$. Since for any $R$-module $M$,
$\Gamma_{\mathcal{Z}}(M)=\bigcup_{\fb\in
F(\mathcal{Z})}\Gamma_{\fb}(M)$, one can easily check that
$H_{\mathcal{Z}}^i(X)\cong {\varinjlim}_{\fb\in
F(\mathcal{Z})}H_{\fb}^i(X)$ for all integers $i$. Hence
$H_{\mathcal{Z}}^i(X)=0$ for all $i>\dim_RX$.

Let $(R,\fm)$ be a local ring, $\mathcal{Z}$ a stable under specialization
subset of
$\Spec R$ and $X\in \mathcal{D}_{\Box}^f(R)$. We prove that
$H_{\mathcal{Z}}^{\dim_RX}(X)=0$ if and only if for any $\fp\in
\Assh_{\hat{R}}(X\otimes_R \hat{R})$, there is $\fq\in \mathcal{Z}$
such that $\dim \hat{R}/\fq \hat{R}+\fp>0$. Yoshino and Yoshizawa \cite[Theorem 2.10]{YY}
have showed that for any abstract local cohomology functor
$\delta:\mathcal{D}_{\sqsubset}(R)\lo \mathcal{D}_{\sqsubset}(R)$,
there is a stable under specialization subset $\mathcal{Z}$ of $\Spec R$
such that $\delta\cong {\bf R}\Gamma_{\mathcal{Z}}$. Thus our result may
be considered as the largest generalization possible of the Hartshorne-Lichtenbaum
Vanishing Theorem. In fact, we show that it includes all known
generalizations of the Hartshorne-Lichtenbaum Vanishing Theorem.

\section{Results}

Let $\mathcal{Z}$ be a stable under specialization subset of $\Spec R$ and $X\in
\mathcal{D}(R)$. The Propositions 2.1 and 2.3 below determine some
situations where the local cohomology modules $H_{\mathcal{Z}}^i(X)$
are Artinian. Recall that $\Supp_RX$ is defined by
$\Supp_RX:=\{\fp\in \Spec R|X_{\fp}\not\simeq 0\}(=\underset{i\in
\mathbb{Z}}\bigcup \Supp_RH_i(X))$.

\begin{proposition} Let $\mathcal{Z}$ be a stable under specialization subset of $\Spec R$ and
$X\in \mathcal{D}_{\sqsubset}^f(R)$. Assume that $\Supp_RX\cap
\mathcal{Z}$ consists only of finitely many maximal ideals. Then
$H_{\mathcal{Z}}^i(X)$ is Artinian for all $i\in \mathbb{Z}$.
\end{proposition}

\begin{prf} Let $\fp$ be a prime ideal and $E(R/\fp)$ denote the
injective envelope of $R/\fp$. Since, $\fp$ is the only associated
prime ideal of $E(R/\fp)$, it turns out
$$\Gamma_{\mathcal{Z}}(E(R/\fp))=\begin{cases} E(R/\fp) & ,\fp\in \mathcal{Z}\\
0& ,\fp\notin \mathcal{Z}.
\end{cases}$$
For each integer $i$, the $R_{\fp}$-module
$\Ext_{R_{\fp}}^i(R_{\fp}/\fp R_{\fp},X_{\fp})$ is finitely
generated, and so $$\mu^i(\fp,X):=\Vdim_{R_{\fp}/\fp
R_{\fp}}\Ext_{R_{\fp}}^i(R_{\fp}/\fp R_{\fp},X_{\fp})<\infty.$$ By
\cite[Proposition 3.18]{F2}, $X$ possesses an injective resolution $I$ such
that $I_i \cong \underset{\fp\in \Spec R}\bigsqcup
E(R/\fp)^{(\mu^i(\fp,X))}$  for all integers $i$. Let $i\in
\mathbb{Z}.$ Then
$$\Gamma_{\mathcal{Z}}(I_i)=\underset{\fp\in \Spec R}\bigsqcup
\Gamma_{\mathcal{Z}}(E(R/\fp))^{(\mu^i(\fp,X))}=\underset{\fp\in
\Supp_RX\cap \mathcal{Z}}\bigsqcup E(R/\fp)^{(\mu^i(\fp,X))}.$$

By the assumption, $\Supp_RX\cap \mathcal{Z}$ consists only of
finitely many maximal ideals. This yields that
$\Gamma_{\mathcal{Z}}(I_i)$ is an Artinian $R$-module, and so
$H_{\mathcal{Z}}^i(X)=H_{-i}(\Gamma_{\mathcal{Z}}(I))$ is Artinian too.
\end{prf}

We record the following immediate corollary which extends
\cite[Theorem 2.2]{Z}. We first recall some definitions. The left
derived tensor product functor $-\otimes_R^{{\bf L}}\sim$ is
computed by taking a projective resolution of the first argument
or of the second one. Also, the right derived homomorphism functor
${\bf R}\Hom_R(-,\sim)$ is computed by taking a projective resolution of
the first argument or by taking an injective resolution of
the second one. Let $\fa$ be an ideal of $R$ and $M,N$
two $R$-modules. The notion of generalized local cohomology modules
$H_{\fa}^i(M,N):={\varinjlim}_n\Ext_R^i(M/\fa^n M,N)$ was introduced
by Herzog in his Habilitationsschrift \cite{He}. When M is finitely
generated, \cite[Theorem 3.4]{Y} yields that $H^i_{\fa}(M,N)\cong
H_{-i}({\bf R}\Gamma_{\fa}({\bf R}\Hom_R(M,N)))$ for all integers
$i$.

\begin{corollary} Let $\fa$ be an ideal of $R$ and $M$ and $N$ two
finitely generated $R$-modules. Assume that $\Supp_RM\cap
\Supp_RN\cap \V(\fa)$ consists only of finitely many maximal ideals.
Then $H_{\fa}^i(M,N)$ is Artinian for all $i\in \mathbb{Z}$.
\end{corollary}

\begin{proposition} Let $\mathcal{Z}$ be a stable under specialization
subset of $\Spec R$. Assume that
for any finitely generated $R$-module $M$ of finite dimension,
$H_{\mathcal{Z}}^{\dim_RM}(M)$ is Artinian. Then for any finite
dimensional complex $X\in \mathcal{D}_{\Box}^f(R)$,
$H_{\mathcal{Z}}^{\dim_RX}(X)$ is Artinian.
\end{proposition}

\begin{prf}  Set $d:=\dim_RX$ and $s:=\sup X$. Clearly, we may assume
that $X\not\simeq 0$, and so $n:=\amp X$ is a non-negative integer.
We argue by induction on $n$. Let $n=0$. Then $X\simeq \Sigma^s
H_s(X)$, and so
$$H_{\mathcal{Z}}^d(X)=H_{\mathcal{Z}}^d(\Sigma^s H_s(X))=H_{\mathcal{Z}}^{d+s}(H_s(X)).$$
On the other hand, by \cite[Proposition 3.5]{F3}, $d=\sup\{\dim H_i(X)-i|i\in
\mathbb{Z}\}$. Hence $\dim H_s(X)=d+s$, and so
$H_{\mathcal{Z}}^{d+s}(H_s(X))$ is Artinian by our assumption. Now,
assume that $n\geq 1$ and let $W:=\tau_{\supset s}X$ and
$Y:=\tau_{s-1\subset}X$ be truncated complexes of $X$; see
\cite[A.1.14]{C}. Since $\amp W=0$ and $\amp Y\leq n-1$, theses
complexes satisfy the induction hypothesis. Next, one has
$$\begin{array}{ll} \dim_RX&=\sup\{\dim_RH_i(X)-i |i\in \mathbb{Z}\}\\
&=\max\{\sup\{\dim_RH_i(X)-i|i\in \mathbb{Z}-\{s\}\} , \dim_RH_s(X)-s\}\\
&=\max\{\dim_RY,\dim_RW\}.
\end{array}$$ Thus by Grothendieck's Vanishing Theorem and induction hypothesis, we deduce
that $H_{\mathcal{Z}}^{d}(W)$ and $H_{\mathcal{Z}}^{d}(Y)$ are
Artinian. Now, by \cite[Theorem 1.41]{F1}, there is a short exact
sequence
$$0\longrightarrow W\longrightarrow X\longrightarrow Y\longrightarrow 0$$
of complexes which induces a long exact sequence
$$H_{\mathcal{Z}}^{d-1}(Y)\longrightarrow H_{\mathcal{Z}}^{d}(W)
\overset{f}\longrightarrow H_{\mathcal{Z}}^{d}(X)\overset{g}
\longrightarrow H_{\mathcal{Z}}^{d}(Y)\longrightarrow 0.$$ It
implies that $H_{\mathcal{Z}}^{d}(X)$ is Artinian.
\end{prf}

\begin{corollary} Let $\mathcal{Z}$ be a stable under specialization subset of $\Spec R$, $\fa$ an
ideal of $R$ and $X\in \mathcal{D}_{\Box}^f(R)$.
\begin{enumerate}
\item[i)]If $R$ is local, then $H_{\mathcal{Z}}^{\dim_RX}(X)$ is Artinian.
\item[ii)]If $\dim_RX$ is finite, then $H_{\fa}^{\dim_RX}(X)$ is Artinian.
\end{enumerate}
\end{corollary}

\begin{prf} In view of the above proposition, i) follows by
\cite[Theorem 2.6 and Lemma 3.2]{DNT} and ii) follows by
\cite[Exercise 7.1.7]{BS}.
\end{prf}

Let $A$ be an Artinian $R$-module. Recall that the set of attached prime ideals of
$A$, $\Att_RA$, is the set of all prime ideals $\fp$ of $R$ such that
$\fp=\Ann_RL$ for some quotient $L$ of $A$. Clearly, $A=0$ if and only if $\Att_RA$ is
empty. If $R$ is local with the maximal ideal $\fm$, then $\Att_RA=\Ass_R(\Hom_R(A,E(R/\fm)))$.
Also, for an exact sequence $0\lo U\lo V\lo W\lo 0$ of Artinian $R$-modules, one can
see $\Att_RW\subseteq \Att_RV\subseteq \Att_RU\cup \Att_RW$. For proving our theorem,
we need to the following lemmas.

\begin{lemma} Let $(R,\fm)$ be a local ring and $X\in \mathcal{D}_{\Box}^f(R)$. Then
$\Att_R(H_{\fm}^{\dim_RX}(X))=\Assh_RX$.
\end{lemma}

\begin{prf} Set $d:=\dim_RX$. By Proposition 2.1, $H_{\fm}^d(X)$ is
an Artinian $R$-module. Hence, we have a natural isomorphism
$H_{\fm}^d(X)\cong H_{\fm}^d(X)\otimes_R\hat{R}$,
and so \cite [Corollary 3.4.4]{L} provides a natural $\hat{R}$-isomorphism
$H_{\fm}^d(X)\cong H_{\fm \hat{R}}^d(X\otimes_R\hat{R})$. From the
definition of attached prime ideals, it follows that
$$\Att_R(H_{\fm}^d(X))=\{\fq\cap R|\fq\in \Att_{\hat{R}}(H_{\fm
\hat{R}}^d(X\otimes_R\hat{R})) \}.$$ Let $\fq$ be a prime ideal of
$\hat{R}$, $\fp:=\fq\cap R$ and $M$ an $R$-module. We have the
natural isomorphism $M_{\fp}\otimes_{R_{\fp}}(\hat{R})_{\fq}\cong
(M\otimes_R\hat{R})_{\fq}$. Since, the natural ring homomorphism
$R_{\fp}\lo (\hat{R})_{\fq}$ is faithfully flat, $M_{\fp}=0$ if and
only if $(M\otimes_R\hat{R})_{\fq}=0$. This implies that $\inf
X_{\fp}=\inf (X\otimes_R\hat{R})_{\fq}$. On the other hand, one can
easily check that $\dim_RX=\dim_{\hat{R}}(X\otimes_R\hat{R})$. Thus,
we can immediately verify that $$\Assh_RX=\{\fq\cap R|\fq\in
\Assh_{\hat{R}}(X\otimes_R\hat{R})\}.$$ Therefore, we may and do
assume that $R$ is complete, and so it possesses a normalized
dualizing complex $D$. By \cite[Chapter V, Theorem 6.2]{Ha}, there
is a natural isomorphism
$$H_{\fm}^i(X)\cong \Hom_R(\Ext_{R}^{-i}(X,D),E(R/\fm))$$ for all integers i.
Since all homology modules of $X$ and of $D$ are finitely generated,
$X$ is homologically bounded and the injective dimension of $D$ is finite, it
follows that  ${\bf R}\Hom_R(X,D)\in \mathcal{D}_{\Box}^f(R)$. In particular,
$\Ext_{R}^{-i}(X,D)$ is a finitely generated $R$-module for all
$i\in \mathbb{Z}$. Thus we have
$$\begin{array}{ll} \Att_R(H_{\fm}^d(X))&=\Att_R(\Hom_R(\Ext_{R}^{-d}(X,D),E(R/\fm)))\\
&=\Ass_R(\Hom_R(\Hom_R(\Ext_{R}^{-d}(X,D),E(R/\fm))),E(R/\fm))\\
&=\Ass_R(\Ext_{R}^{-d}(X,D))\\
&=\Ass_R(H_d({\bf R}\Hom_R(X,D))).
\end{array}
$$
\cite [Theorem 16.20]{F1} implies that $\sup ({\bf
R}\Hom_R(X,D))=d$. Let $\fp\in \Spec R$. By \cite [Theorem
12.26]{F1}, $\fp \in \Ass_R(H_d({\bf R}\Hom_R(X,D)))$ if and only if
$\depth_{R_{\fp}}{\bf R}\Hom_R(X,D)_{\fp}=-d$. But, \cite[Lemma A.6.4 and A.6.32]{C}
and \cite [Theorem 15.17]{F1}, yield that

$$\begin{array}{ll} \depth_{R_{\fp}}{\bf R}\Hom_R
(X,D)_{\fp}&=\depth_{R_{\fp}}({\bf
R}\Hom_{R_{\fp}}(X_{\fp},D_{\fp}))\\
&=\depth_{R_{\fp}}D_{\fp}+ \inf X_{\fp}\\
&=-\dim \frac{R}{\fp}+\inf X_{\fp}.
\end{array}$$
Therefore,  $\fp\in \Ass_R(H_d({\bf R}\Hom_R(X,D)))$ if and only if
$\dim\frac{R}{\fp}-\inf X_{\fp}=\dim_RX$. This means
$\Att_R(H_{\fm}^d(X))=\Assh_RX$, as desired.
\end{prf}

\begin{lemma} Let $(R,\fm)$ be a local ring, $\mathcal{Z}$
a stable under specialization subset of $\Spec R$ and $X\in \mathcal{D}_{\Box}^f(R)$.
Then $H_{\mathcal{Z}}^{\dim_RX}(X)$ is a homomorphic image of
$H_{\fm}^{\dim_RX}(X)$.
\end{lemma}

\begin{prf} Let $\fa$ be an ideal of $R$ and $x\in \fm$.
Let $I$ be an injective resolution of $X$. Then $I_x$, the
localization of $I$ at $x$, provides an injective resolution of
$X_x$ in $\mathcal{D}_{\Box}^f(R_x)$. Now, \cite[Lemma 8.1.1]{BS}
yields the following  exact sequence of complexes
$$0\lo \Gamma_{\fa+(x)}(I)\lo \Gamma_{\fa}(I)\lo \Gamma_{\fa}(I_x)\lo 0,$$
where the maps are the natural ones. Set $d:=\dim_RX$. We deduce the
long exact sequence
$$\cdots \lo H_{\fa+(x)}^d(X)\lo H_{\fa}^d(X)\lo
H_{\fa R_x}^d(X_x)\lo 0.$$ By Corollary 2.4, $H_{\fa}^d(X)$ is
Artinian. Hence $H_{\fa}^d(X)$ is supported at most at $\fm$, and so
$$H_{\fa R_{x}}^d(X_{x})\cong H_{\fa}^d(X) _{x}=0.$$ Hence, the
natural homomorphism $H_{\fa+(x)}^d(X)\lo H_{\fa}^d(X)$ is epic.

We may choose  $x_1,x_2, \dots ,x_n \in R$ such that
$\fm=\fa+(x_1,x_2,\dots ,x_n)$. Set $\fa_i: =\fa+(x_1,\dots
,x_{i-1})$ for $i=1,\dots , n+1$. By the above argument, the natural
homomorphism $H_{\fa_{i+1}}^d(X)\lo H_{\fa_{i}}^d(X)$ is epic for
all $1\leq i\leq n$. Hence $H_{\fa}^d(X)$ is a homomorphic image of
$H_{\fm}^d(X)$. This completes the proof,  because
$H_{\mathcal{Z}}^d(X)\cong \underset{\fb}{\varinjlim}H_{\fb}^d(X)$,
where the direct limit is over all ideals $\fb$ of $R$ such that
$\V(\fb)\subseteq \mathcal{Z}$.
\end{prf}

\begin{lemma} Let $M$ be a finitely generated $R$-module and $X\in \mathcal{D}_{\Box}^f(R)$.
\begin{enumerate}
\item[i)] $\dim_R(M\otimes_R^{{\bf L}}X)\leq \dim_RX$.
\item[ii)] If $\Supp_RM\cap \Assh_RX\neq \emptyset$, then $\dim_R(M\otimes_R^{{\bf L}}X)=\dim_RX$ and
$$\Assh_R(M\otimes_R^{{\bf L}}X)=\Supp_RM\cap \Assh_RX.$$
\end{enumerate}
\end{lemma}

\begin{prf} For any Noetherian local ring S and any two complexes
$V,W\in \mathcal{D}_{\Box}^f(S)$, Nakayama's Lemma for complexes
asserts that $\inf(V\otimes_R^{{\bf L}}W)=\inf V+\inf W$; see e.g.
\cite[Corollary A.4.16]{C}. In particular, this yields that
$\Supp_R(V\otimes_R^{{\bf L}}W)=\Supp_RV\cap \Supp_RW$. Now, by noting that
for any complex $Y\in \mathcal{D}(R)$, we have $$\dim_RY=\sup \{\dim
R/\fp-\inf Y_{\fp}|\fp\in \Supp_RY\},$$ both assertions follow
immediately.
\end{prf}

Next, we conclude our theorem.

\begin{theorem} Let $(R,\fm)$ be a local ring, $\mathcal{Z}$ a stable under specialization
subset of $\Spec R$ and $X\in \mathcal{D}_{\Box}^f(R)$. Then
$\Att_{\hat{R}}(H_{\mathcal{Z}}^{\dim_RX}(X))= \{\fp\in
\Assh_{\hat{R}}(X\otimes_R \hat{R})|\dim \hat{R}/\fq \hat{R}+\fp=0 \
\ \text{for all} \  \ \fq\in \mathcal{Z}\}$.
\end{theorem}

\begin{prf} Set $d:= \dim_RX$ and $s:=\sup X$. We may assume that
$n:=\amp X$ is a non-negative integer. First, by induction on $n$,
we prove the inclusion $\subseteq$. If $n=0$, then $X\simeq \Sigma ^s
H_s(X)$, and so $$H_{\mathcal{Z}}^d(X)=H_{\mathcal{Z}}^d(\Sigma ^s
H_s(X))=H_{\mathcal{Z}}^{d+s}(H_s(X)).$$ In the proof of Proposition
2.3, we saw that $\dim H_s(X)= d+s$, hence \cite [Corollary
2.7]{DNT} implies that

$$\begin{array}{ll} \Att_{\hat{R}}(H_{\mathcal{Z}}^d(X))&=
\Att_{\hat{R}}(H_{\mathcal{Z}}^{d+s}(H_s(X)))\\
&=\{\fp\in \Assh_{\hat{R}}(H_s(X)\otimes_R\hat{R})|\dim \hat{R}/\fq
\hat{R}+\fp=0 \ \ \text{for all} \  \
\fq\in \mathcal{Z}\}\\
&=\{\fp\in \Assh_{\hat{R}}(X\otimes_R \hat{R})|\dim \hat{R}/\fq
\hat{R}+\fp=0 \ \ \text{for all} \  \ \fq\in \mathcal{Z}\}.
\end{array}$$
Now, assume that $n\geq 1$ and $\fp\in
\Att_{\hat{R}}(H_{\mathcal{Z}}^d(X))$. By Lemma 2.6,
$H_{\mathcal{Z}}^d(X)$ is an homomorphic image of $H_{\fm}^d(X)$,
and so Lemma 2.5 yields that
$$\Att_{\hat{R}}(H_{\mathcal{Z}}^d(X))\subseteq
\Att_{\hat{R}}(H_{\fm}^d(X))=\Assh_{\hat{R}}(X\otimes_R \hat{R}).$$
Thus $\fp \in\Assh_{\hat{R}}(X\otimes_R \hat{R})$. Let
$W:=\tau_{\supset s}X$ and $Y:=\tau_{s-1\subset}X$ be truncated
complexes of $X$. We have a short exact sequence
$$0\lo W\lo X\lo Y\lo 0$$
of complexes and from the proof of Proposition 2.3, we know that
$\dim_RX=\max\{\dim_RW,\dim_RY\}$. From the long exact sequence
$$\cdots \lo H_{\mathcal{Z}}^d(W)\lo
H_{\mathcal{Z}}^d(X)\lo H_{\mathcal{Z}}^d(Y)\lo 0,$$ we deduce that
$$\Att_{\hat{R}}(H_{\mathcal{Z}}^d(X))\subseteq
\Att_{\hat{R}}(H_{\mathcal{Z}}^d(W))\cup
\Att_{\hat{R}}(H_{\mathcal{Z}}^d(Y)).$$ Thus, either
$\fp\in\Att_{\hat{R}}(H_{\mathcal{Z}}^d(W))$ or
$\fp\in\Att_{\hat{R}}(H_{\mathcal{Z}}^d(Y))$. By Grothendieck's
Vanishing Theorem, the first case implies that $\dim_RW=d$ and the
second case implies that $\dim_RY=d$. Since $\amp W=0$ and $\amp
Y\leq n-1$, in both cases, the induction hypothesis  yields that
$\dim \hat{R}/\fq \hat{R}+\fp=0$ for all $\fq\in \mathcal{Z}$.

Now, we prove the inclusion $\supseteq$. Let $\fp
\in\Assh_{\hat{R}}(X\otimes_R \hat{R})$ be such that $\dim
\hat{R}/\fq \hat{R}+\fp=0 \ \ \text{for all} \ \ \fq\in
\mathcal{Z}$. We have to show that
$\fp\in\Att_{\hat{R}}(H_{\mathcal{Z}}^d(X))$. Since
$H_{\mathcal{Z}}^d(X)$ is an Artinian $R$-module, we have the
natural isomorphism $H_{\mathcal{Z}}^d(X)\cong
H_{\mathcal{Z}}^d(X)\otimes_R\hat{R}$. On the other hand, by \cite
[Corollary 3.4.4] {L}, for any ideal $\fa$ of $R$, there is a
natural $\hat{R}$-isomorphism $H_{\fa}^d(X)\otimes_R\hat{R}\cong
H_{\fa \hat{R}}^d(X\otimes_R\hat{R})$. Let $\hat{\mathcal{Z}}:=
\{\fq\in \Spec \hat{R}|\fq\cap R\in
\mathcal{Z}\}$, which can be easily checked that is
a stable under specialization subset of $\Spec \hat{R}$. It is
straightforward to see that the two families
$\{\fa \hat{R}|\fa \  \ \text{is an ideal of} \  \ R \  \  \text{with} \  \
\V(\fa)\subseteq \mathcal{Z}\}$ and $\{\fb|\fb \  \ \text{is an ideal of}
\  \ \hat{R} \   \ \text{with} \  \
\V(\fb)\subseteq \hat{\mathcal{Z}}\}$ are cofinal. This implies
that $H_{\mathcal{Z}}^d(X)\cong H_{\hat{\mathcal{Z}}}^d(X\otimes_R\hat{R})$.
Also, we have $\dim_{\hat{R}}(X\otimes_R\hat{R})=\dim_RX$ and $\dim
\hat{R}/\fq+\fp=0$ for all $\fq\in \hat{\mathcal{Z}}$.
Therefore,  we may and do assume that $R$ is complete.

Since $R$ is complete, there is a complete regular local ring $(T,\fn)$ and a surjective
ring homomorphism $f:T\lo R$. One can easily check that $X\in \mathcal{D}_{\Box}^f(T)$ and
$\dim_TX=\dim_RX$. Set $\bar{\mathcal{Z}}:=\{f^{-1}(\fq)|\fq\in \mathcal{Z}\}$,
which is clearly a stable under specialization subset of $\Spec T$. By \cite
[Corollary 3.4.3] {L}, for any ideal $\fb$ of $T$, there is a natural
$T$-isomorphism $H_{\fb}^d(X)\cong H_{\fb R}^d(X)$. From this, we can conclude a
natural $T$-isomorphism $H_{\bar{\mathcal{Z}}}^d(X)\cong H_{\mathcal{Z}}^d(X)$.
For any Artinian $R$-module $A$ and any $\fq\in \Spec R$, it turns out that $A$
is also Artinian as a $T$-module and $\fq\in \Att_RA$ if and only if $f^{-1}(\fq)\in
\Att_TA$. Finally, we have $\dim T/\bar{\fq}+f^{-1}(\fp)=0$ for all $\bar{\fq}\in \bar{\mathcal{Z}}$
and $\Assh_TX=\{f^{-1}(\fq)|\fq\in \Assh_RX\}$.
Thus from now on, we can assume that $R$ is a complete regular local ring.

Lemma 2.7 yields that $\dim_R(\fp\otimes_R^{{\bf L}}X)\leq\dim_RX$ and
$\dim_R(R/\fp\otimes_R^{{\bf L}}X)=\dim_RX$. Let $P$ be a projective resolution of $X$.
Applying $-\otimes_RP$ to the short exact
sequence $$0\lo \fp\lo R \lo R/\fp\lo 0,$$ yields the following
exact sequence of complexes
$$0\lo \fp\otimes_R^{{\bf L}}X\lo X\lo R/\fp\otimes_R^{{\bf L}}X\lo 0.$$ It yields
the following exact sequence
$$\cdots \lo H_{\mathcal{Z}}^d(\fp\otimes_R^{{\bf L}}
X)\lo H_{\mathcal{Z}}^d(X )\lo H_{\mathcal{Z}}^d(R/\fp\otimes_R^{{\bf L}}
X)\lo 0.$$ As $R$ is regular, the projective dimension of any $R$-module is finite, and so for any
finitely generated $R$-module $M$, one has $M\otimes_R^{{\bf L}}X\in \mathcal{D}_{\Box}^f(R)$.
Since $\dim R/\fq+\fp=0$ for all $\fq\in \mathcal{Z}$, it
follows that $\Gamma_{\mathcal{Z}}(\Gamma_{\fp}(M))=\Gamma_{\fm}(M)$
for all $R$-modules $M$. Let $I$ be an injective resolution of
$R/\fp\otimes_R^{{\bf L}} X$. Since
$$\Supp_RI=\Supp_R(R/\fp\otimes_R^{{\bf L}}X)\subseteq \V(\fp),$$ by
\cite [Corollary 3.2.1]{L}, $\Gamma_{\fp}(I)\simeq I$, and so
$$\Gamma_{\mathcal{Z}}(I)\simeq \Gamma_{\mathcal{Z}}(\Gamma_{\fp}(I))=\Gamma_{\fm}(I).$$
In particular, there is an isomorphism
$H_{\mathcal{Z}}^d(R/\fp\otimes_R^{{\bf L}} X)\cong
H_{\fm}^d(R/\fp\otimes_R^{{\bf L}} X)$. Therefore, by Lemmas 2.7 and 2.5, we
deduce that $\fp\in\Att_R(H_{\mathcal{Z}}^d(R/\fp\otimes_R^{{\bf L}} X))
\subseteq \Att_R(H_{\mathcal{Z}}^d(X))$.
\end{prf}

Now, we are ready to establish the derived category analogue  of the
Hartshorne-Lichtenbaum Vanishing Theorem.

\begin{corollary} Let $(R,\fm)$ be a local ring, $\mathcal{Z}$ a stable under specialization
subset of $\Spec R$ and $X\in \mathcal{D}_{\Box}^f(R)$.
The following are equivalent:
\begin{enumerate}
\item[i)] $H_{\mathcal{Z}}^{\dim_RX}(X)=0$.
\item[ii)] For any $\fp\in \Assh_{\hat{R}}(X\otimes_R \hat{R})$, there is
$\fq\in \mathcal{Z}$ such that $\dim \hat{R}/\fq \hat{R}+\fp>0$.
\end{enumerate}
\end{corollary}

\begin{corollary} Let $\fa$ be an ideal of the local ring $(R,\fm)$ and $X\in
\mathcal{D}_{\Box}^f(R)$.\\
1)  $\Att_{\hat{R}}(H_{\fa}^{\dim_RX}(X))= \{\fp\in
\Assh_{\hat{R}}(X\otimes_R \hat{R})|\dim
\hat{R}/\fa \hat{R}+\fp=0 \}$.\\
2) The following are equivalent:
\begin{enumerate}
\item[i)] $H_{\fa}^{\dim_RX}(X)=0$.
\item[ii)] $\dim \hat{R}/\fa \hat{R}+\fp>0$ for all $\fp\in
\Assh_{\hat{R}}(X\otimes_R \hat{R})$.
\end{enumerate}
\end{corollary}

\begin{corollary} Let $(R,\fm)$ be a local ring, $\mathcal{Z}$ a stable under specialization
subset of $\Spec R$ and $M,N$ two finitely generated $R$-modules.
Assume that ${\bf R}\Hom_R(M,N)\in \mathcal{D}_{\Box}^f(R)$ and set
$d:=\dim_R({\bf R}\Hom_R(M,N))$. The following are equivalent:
\begin{enumerate}
\item[i)] $H_{\mathcal{Z}}^d(M,N)=0$.
\item[ii)] For any $\fp\in \Assh_{\hat{R}}({\bf R}\Hom_{\hat{R}}(\hat{M},\hat{N}))$, there is
$\fq\in \mathcal{Z}$ such that $\dim \hat{R}/\fq \hat{R}+\fp>0$.
\end{enumerate}
\end{corollary}

\begin{prf} Note that ${\bf R}\Hom_R(M,N)\otimes_R\hat{R}\simeq {\bf R}\Hom_{\hat{R}}(\hat{M},
\hat{N})$, and so the result follows by Corollary 2.9.
\end{prf}

\begin{remark} Let $\mathcal{Z}$ be a stable under specialization subset of $\Spec R$ and
 $X\in \mathcal{D}_{\Box}^f(R)$.
\begin{enumerate}
\item[1)] Suppose that dimension of $X$ is finite.
Then $H_{\mathcal{Z}}^{\dim_RX}(X)$ is supported only at maximal
ideals. To realize this, in view of Corollary 2.4 ii), it is enough
to notice that $H_{\mathcal{Z}}^{\dim_RX}(X)\cong
{\varinjlim}_{\fa}H_{\fa}^{\dim_RX}(X)$, where the direct limit is
over all ideals $\fa$ of $R$ such that $\V(\fa)\subseteq
\mathcal{Z}$. But, $H_{\mathcal{Z}}^{\dim_RX}(X)$ is not Artinian in
general.  To this end, let $R$ be a finite dimensional Gorenstein
ring such that the set $\mathcal{Z}:=\{\fm\in \Max R|\Ht \fm=\dim R\}$ is
infinite. Clearly, $\mathcal{Z}$ is a stable under specialization subset of $\Spec R$.
The minimal injective resolution of $R$ has the form
$$0\lo \bigsqcup_{\Ht \fp=0}
E(R/\fp)\lo \bigsqcup_{\Ht \fp=1} E(R/\fp)\lo \cdots \lo
\bigsqcup_{\Ht \fp=\dim R} E(R/\fp)\lo 0.$$ Hence
$H_{\mathcal{Z}}^{\dim R}(R)= \bigsqcup_{\fm\in\mathcal{Z}}
E(R/\fm)$, which is not Artinian.
\item[2)] Suppose that $R$ is local with the maximal ideal $\fm$
and $I,J$ two ideals of $R$. In \cite{TYY}, Takahashi, Yoshino and Yoshizawa considered the
following stable under specialization subset of $\Spec R$
$$W(I,J)=\{\fp\in \Spec(R)|I^n\subseteq \fp+J \   \ \text{for a natural integer n}\}.$$
For each integer $i$, they called $H_{I,J}^i(-):=H^i_{W(I,J)}(-)$,
$i$-th local cohomology functor with respect to $(I,J)$. For the
ring $R$ itself, they extended the Hartshorne-Lichtenbaum Vanishing
Theorem ; see \cite[Theorem 4.9]{TYY}. Namely, they showed that
$H^{\dim R}_{I,J}(R)=0$ if and only if for any prime ideal $\fp\in
\Assh_{\hat{R}}\hat{R}\cap \V(J \hat{R})$, we have $\dim \hat{R}/I
\hat{R}+\fp>0$. On the other hand by \cite[Theorem 2.8]{DNT},
$H^{\dim R}_{I,J}(R)=0$ if and only if for any  prime ideal $\fp\in
\Assh_{\hat{R}}\hat{R}$, there is $\fq\in W(I,J)$ such that $\dim
\hat{R}/\fq \hat{R}+\fp>0$. Hence the following statements are
equivalent:\\
i) For any prime ideal $\fp\in \Assh_{\hat{R}}\hat{R}\cap \V(J
\hat{R})$, we have $\dim \hat{R}/I \hat{R}+\fp>0$. \\
ii) For any prime ideal $\fp\in \Assh_{\hat{R}}\hat{R}$, there is
$\fq\in W(I,J)$ such that $\dim \hat{R}/\fq \hat{R}+\fp>0$.\\
As  Takahashi, Yoshino and Yoshizawa  \cite[Remark 4.10]{TYY} have
mentioned, it is not so easy to check the equivalence of these
statements directly. Here, we do this under the extra assumption that $R$ is
complete. (In fact this assumption is not needed for the implication
$ii)\Longrightarrow i)$.) Suppose $ii)$ holds and let
$\fp\in \Assh_{R}R\cap \V(J)$. By the assumption
there is $\fq\in W(I,J)$ such that $\dim R/\fq+\fp>0$.
Since $\fq\in W(I,J)$, there is a natural integer $n$, such that
$I^n\subseteq \fq+J$. This yields that $I^n+\fp \subseteq
\fq+\fp,$ and so $$\dim R/I+\fp=\dim R/I^n+\fp\geq \dim R/\fq+\fp>0.$$
Conversely, suppose that $i)$ holds and let $\fp\in
\Assh_{R}R$. First, assume that $J\subseteq
\fp$. Then by the assumption, $\dim R/I+\fp>0$, and so
there is $\fq\in \V(I+\fp)$ such that $\dim R/\fq>0$.
Then $I+J\subseteq I+\fp\subseteq \fq$. Hence $\fq\in W(I,J)$ and $\dim
R/\fq+\fp=\dim R/\fq>0$. Thus $ii)$
follows when $J\subseteq \fp$. Now, assume that
$J\nsubseteq \fp$. By \cite[Lemma 3.3]{TYY}, $$\V(J)=
\bigcap_{\fq\in W(\fm,J)}W(\fm,\fq).$$ So, there is $\fq\in W(\fm,J)
\subseteq W(I,J)$ such that $\fp\notin W(\fm,\fq)$. Since $\fp\notin W(\fm,\fq)$,
it follows that $\fq+\fp$ is not $\fm$-primary, and so $\dim R/\fq+\fp>0$.
\item[3)] Suppose that $R$ is local and  $F(\mathcal{Z})$ denote
the set of all ideals $\fb$ of $R$ such that $\V(\fb)\subseteq
\mathcal{Z}$. As we mentioned in the introduction
$H_{\mathcal{Z}}^i(X)\cong {\varinjlim}_{\fb\in
F(\mathcal{Z})}H_{\fb}^i(X)$ for all integers $i$. The relationship
between $H_{\mathcal{Z}}^{\dim_RX}(X)$ and $H_{\fb}^{\dim_RX}(X)'$s
is more deeper. In fact by Theorem 2.8, we have
$$\Att_{\hat{R}}(H_{\mathcal{Z}}^{\dim_RX}(X))=\bigcap_{\fb\in
F(\mathcal{Z})}\Att_{\hat{R}}(H_{\fb}^{\dim_RX}(X)).$$ This implies
that $H_{\mathcal{Z}}^{\dim_RX}(X)=0$ if and only if
$H_{\fb}^{\dim_RX}(X)=0$ for an ideal $\fb\in F(\mathcal{Z})$.
\end{enumerate}
\end{remark}

\begin{acknowledgement}Part of this research was done during the second
author's visit to the Department of Mathematics at the University of
Nebraska-Lincoln. He thanks this department for its kind
hospitality.
\end{acknowledgement}


\end{document}